\documentclass[12pt]{amsart}

\theoremstyle{plain}
\newtheorem{theorem}{Theorem}
\newtheorem{corollary}[theorem]{Corollary}
\newtheorem{lemma}[theorem]{Lemma}
\newtheorem{proposition}[theorem]{Proposition}

\theoremstyle{definition}
\newtheorem{definition}[theorem]{Definition}
\newtheorem{remark}[theorem]{Remark}

\newtheorem{example}[theorem]{Example}

\newcommand{\abs}[1]{\lvert#1\rvert}
\newcommand{\norm}[1]{\lVert#1\rVert}
\newcommand{\bigabs}[1]{\bigl\lvert#1\bigr\rvert}
\newcommand{\bignorm}[1]{\bigl\lVert#1\bigr\rVert}

\renewcommand{\leq}{\leqslant}
\renewcommand{\geq}{\geqslant}\usepackage{amssymb}
\newcommand{\marg}[1]{}
\newcommand{\lab}[1]{\label{#1}}
\newcommand{\term}[1]{{\textit{\textbf{#1}}}}
\newcommand{\vr}{\varepsilon}
\newcommand{\one}{1}
\newcommand{\raw}{\to}

\def\mid{\::\:}

\def\iS{\mathcal S}
\def\iA{\mathcal A}

\DeclareMathOperator{\Range}{Range}

\begin{document}

\baselineskip=18pt

\title{A theorem of Krein revisited}

\author[T.Oikhberg and V.G.Troitsky]{Timur Oikhberg and Vladimir G. Troitsky}
\address{
 Department of Mathematics, University of California, Irvine, CA
 92697-3875. USA.}
\address{Department of Mathematical and Statistical Sciences,
  University of Alberta, Edmonton, AB, T6G\,2G1. Canada.}
\email{toikhber@math.uci.edu\\vtroitsky@math.ualberta.ca}

\thanks{The research of the first author was supported by
         NSF grant DMS-9970369.}

%\date{\today.}

\begin{abstract}
  M.~Krein proved in~\cite{Krein:48} that if $T$ is a
  continuous operator on a normed space leaving invariant an open
  cone, then its adjoint $T^*$ has an eigenvector. We present
  generalizations of this result as well as some applications to
  $C^*$-algebras, operators on $\ell_1$, operators with invariant
  sets, contractions on Banach lattices, the Invariant Subspace
  Problem, and von Neumann algebras.
\end{abstract}

\maketitle

M.~Krein proved in~\cite[Theorem~3.3]{Krein:48}
that if $T$ is a continuous operator on a normed space leaving
invariant a non-empty open cone, then its adjoint $T^*$ has an eigenvector.
Krein's result has an immediate application to the Invariant Subspace
Problem because of the following observation. If $T$ is a bounded
operator on a Banach space and not a multiple of the identity, and
$T^*f=\lambda f$, then the kernel of $f$ is a closed non-trivial subspace
of codimension 1 which is invariant under~$T$. Moreover,
$\overline{\Range(\lambda I-T)}$ is a closed nontrivial subspace which is
proper (it is contained in the kernel of $f$) and hyperinvariant for
$T$, that is, it is invariant under every operator commuting with
$T$.

Several proofs and modifications of Krein's theorem appear in the
literature, see, e.g.,~\cite[Theorems 6.3 and 7.1]{Abramovich:92}
and~\cite[p.~315]{Schaefer}. We prove yet another version of Krein's
Theorem: if $T$ is a positive operator on an ordered normed space in
which the unit ball has a dominating point, then $T^*$ has a positive
eigenvector. We deduce the original Krein's version of the theorem
from this, as well as several applications and related results.
In particular, we show that if a bounded operator $T$ on a Banach
space satisfies any of the following conditions, then $T^*$ has an
eigenvector. Moreover, if the condition holds for a commutative family
of operators, then the family of the adjoint operators has a common
eigenvector.
\begin{itemize}
  \item $T$ leaves invariant a cone with an interior point;
  \item $T$ is a positive operator on a $C^*$-algebra;
  \item $T$ is an operator on $\ell_1$ such that entries of its matrix
    satisfy $t_{kk}\pm t_{kj}\geq\sum_{i \neq k}\abs{t_{ik} \pm
      t_{ij}}$;
  \item $T$ leaves invariant a convex set whose interior is non-void
    and doesn't contain zero;
  \item $T$ is a contraction with a fixed point;
  \item $T$ is a positive contraction on a Banach lattice and $Te>e$
    for some $e>0$.
\end{itemize}
We also show that under the last condition $T$ has a closed invariant
order ideal. Finally, we prove a non-commutative version of this
result for rearrangement invariant operator spaces arising from von
Neumann algebras.

Throughout the paper $X$ denotes a real or complex normed space, $X^*$
the dual of~$X$, $T$ a bounded linear operator on~$X$, and $B_X$ the
closed unit ball of~$X$.

\begin{definition}
  We call a subset $K$ of a normed space $X$ a \term{cone} if
  $K$ is closed under addition and non-negative scalar multiplication,
  and there exists a non-zero vector $x\in K$ such that $-x\notin K$.
\end{definition}

Our definition of a cone is a most general one. In the literature such
objects are sometimes called \term{wedges}, while for a cone it is
often assumed in addition that $x\in K$ implies $-x\notin K$ for
{\em every} non-zero $x$. This additional condition ensures that the
relation on $X$ defined via ``$x\leq y$ iff $y-x\in K$'' is a linear
order relation, and, vice versa, every linear order relation defines a
cone satisfying this condition, namely, the cone $X_+$ of all
non-negative elements. We will still use the symbol ``$\leq$'', even
though in our case $x>0$ and $x<0$ may happen simultaneously. However,
this does not create any problems, and, naturally, everything we do is
still valid for the more restrictive definitions of a cone.
See~\cite{Krein:48} for a discussion on definitions and
properties of cones.

Given a closed cone $K$ in a normed space $X$, we will call $X$ an
\term{ordered normed space} with respect to the (semi)order relation
determined by~$K$. Notice that $K$ coincides with the cone $X_+$ of
all non-negative elements of $X$. A linear operator is said to be
\term{positive} if $T(X_+)\subseteq X_+$. For $f\in X^*$ we write $f\geq 0$ or
$f\in X^*_+$ if $f(x)\geq 0$ whenever $x\geq 0$. Clearly, $X^*_+$ is a
$w^*$-closed cone in $X^*$.  It can be easily verified that if $T$ is
a positive bounded operator on $X$ then $T^*$ is a positive operator
on $X^*$, that is, $T^*(X^*_+)\subseteq X^*_+$.  It is known (see,
e.g.~\cite{Krein:48}) that if $K$ is a closed cone, then $K$
and $-K$ can be (non-strictly) separated by a continuous functional,
or, equivalently, there exists a non-trivial positive functional in
$X^*$.

\begin{lemma} \lab{l:posfunc}
  Suppose that $X$ is a real normed space and $e\in X$ with
  $\norm{e}=1$. If $f\in X^*$ then $f(e)=\norm{f}$ if and only if
  $f(x)\leq f(e)$ for all $x\in B_X$.
\end{lemma}

\begin{proof}
  If $f(e)=\norm{f}$ then
  $f(x)\leq\abs{f(x)}\leq\norm{f}\norm{x}=f(e)$
  whenever $\norm{x}=1$. Conversely, suppose
  $f(x)\leq f(e)$ for all $x$ of norm one.
  Since $-f(x)=f(-x)\leq f(e)$, we have
  $\abs{f(x)}\leq f(e)$, so that $\norm{f}\leq f(e)$.
  Finally, $f(e)=\abs{f(e)}\leq\norm{f}$.
\end{proof}

\begin{definition}
  If $X$ is an ordered normed space and $e\in B_X$, we say that
  \term{$e$ dominates the unit ball of $X$} if $x\leq e$ for all $x\in
  B_X$. We write $B_X\leq e$ then.
\end{definition}

In this case it follows immediately from Lemma~\ref{l:posfunc} that
every positive functional attains its norm on $e$. In the proof of the
following theorem we use techniques developed in the proof of a
special case of Krein's theorem in \cite{Abramovich:92,AA:02}.

\begin{theorem} \label{t:dom-cone}
  Suppose that $X$ is an ordered real normed space and $e\in X$ such
  that $\norm{e}=1$ and $B_X\leq e$. If $T$ is a positive operator on
  $X$ then $T^*$ has a positive eigenvector. Moreover, if $\Gamma$ is a
  commutative family of positive operators on $X$, then their adjoints
  have a common positive eigenvector.
\end{theorem}

\begin{proof}
  Let $\iS=\{f\in X^*_+\mid f(e)=1\}$. Since $\iS=X^*_+\cap\{f\in X^*\mid
  f(e)=1\}$ then $\iS$ is $w^*$-closed. Furthermore, if $f\in\iS$ then
  $\norm{f}=f(e)=1$ by Lemma~\ref{l:posfunc}, so that $\iS\subseteq B_X$, hence
  is $w^*$-compact.  For $T\geq 0$ and $f\in\iS$ we define
   \begin{displaymath}
      F_T(f)=\frac{f+T^*f}{[f+T^*f](e)}=\frac{f+T^*f}{1+(T^*f)(e)}.
  \end{displaymath}
  Since $T^*\geq 0$ then $F_T(f)\geq 0$. Clearly, $\bigl[F_T(f)\bigr](e)=1$,
  so that $F_T(f)\in\iS$, hence $F_T\colon \iS\to \iS$.
  It can be easily verified that $F_T$ is $w^*$-to-$w^*$-continuous.
  Indeed, if $f_\alpha\xrightarrow{w^*}f$, then for all $x\in X$ we
  have
  \begin{displaymath}
    \left[F_T(f_\alpha)\right](x)=
    \frac{f_\alpha(x)+(T^*f_\alpha)(x)}{1+(T^*f_\alpha)(e)}\to
    \frac{f(x)+(T^*f)(x)}{1+(T^*f)(e)}=\left[F_T(f)\right](x)
  \end{displaymath}
  because $T^*$ is $w^*$-to-$w^*$-continuous.
  By the Fixed Point
  Theorem there exists $h\in \iS$ such that $F_T(h)=h$, i.e.,
  $\frac{h+T^*h}{1+(T^*h)(e)}=h$ so that
  $T^*h=\bigl((T^*h)(e)\bigr)h$, hence $h$ is an eigenvector of $T^*$.
  
  Let $\Gamma$ be a commutative family of positive operators on $X$.
  For $T\in\Gamma$ denote $A_T$ the set of the fixed points of $F_T$
  in $\iS$.  It can be easily verified that $f\in\iS$ belongs to $A_T$
  if and only if $f$ is an eigenvector of $T^*$.  Clearly, $A_T$ is
  $w^*$-closed, hence $w^*$-compact.  We claim that
  $\{A_T\}_{T\in\Gamma}$ has the finite intersection property, this
  would imply that it has non-empty intersection, and, therefore, the
  family $\{T^*\}_{T\in\Gamma}$ has a common eigenvector in $\iS$. We
  prove the claim by induction on the size of the set. Suppose that
  $\bigcap_{T\in\Gamma_0}A_T\neq\varnothing$ for every $n$-element
  subset $\Gamma_0\subseteq\Gamma$. Let $\Gamma_0$ be an $n$-element
  subset of $\Gamma$ and $S\in\Gamma$, show that
  $\bigcap_{T\in\Gamma_0\cup\{S\}}A_T\neq\varnothing$.  Pick
  $f\in\bigcap_{T\in\Gamma_0}A_T$, then for each $T$ in $\Gamma_0$
  there exists $\lambda_T\geq 0$ such that $T^*f=\lambda_Tf$. Let
  $C_T=\ker(\lambda_TI-T^*)\cap\iS$, then $C_T$ is a convex
  $w^*$-closed subset of $A_T$, hence $w^*$-compact. It follows that
  $C=\bigcap_{T\in\Gamma_0}C_T$ is convex and $w^*$-compact.
  Furthermore, $C\neq\varnothing$ as $f\in C$. If $T\in\Gamma_0$ and
  $h\in C_T$, then
  $$T^*F_S(h)=\frac{T^*h+T^*S^*h}{1+(S^*h)(e)}=
  \frac{\lambda_Th+S^*(\lambda_T h)}{1+(S^*h)(e)}=\lambda_TF_S(h),$$
  so that $F_S(h)\in C_T$. It follows that $F_S(C_T)\subseteq C_T$, so
  that $F_S(C)\subseteq C$. The Fixed Point Theorem implies that $F_S$
  has a fixed point in $C$, hence $A_S\cap C\neq\varnothing$. Since
  $C\subseteq\bigcap_{T\in\Gamma_0}A_T$, this proves the claim.
\end{proof}

\begin{theorem} \label{t:Krein}
  If $T$ is a continuous operator on a real normed space, leaving
  invariant a cone with an interior point, then $T^*$ has a positive
  eigenvector. Moreover, a commutative collection of such operators
  has a common positive eigenvector.
\end{theorem}

\begin{proof}
  Let $\Gamma$ be a commutative family of bounded operators on $X$,
  $C$ a cone in $X$ such that $T(C)\subseteq C$ for each $T\in\Gamma$, and
  $e$ an interior point of $C$. Without loss of generality,
  $C$ is closed, $\norm{e}>1$, and $e+B_X\subseteq C$. Let $C_0$ be the cone spanned
  by $e+B_X$, that is, $C_0=\{\alpha(e+x)\mid \alpha\geq
  0,\,\norm{x}\leq 1\}$. Put $W=(C_0-e)\cap(e-C_0)$.
  Note that $e+B_X\subseteq C_0$ so that $B_X\subset C_0-e$.
  Also, $B_X=-B_X\subseteq e-C_0$, so that
  $B_X\subseteq W$. Furthermore, $W$ is bounded. Indeed, if $w\in W$
  then $w=\alpha_1(e+x_1)-e=e-\alpha_2(e+x_2)$ for some
  $\alpha_1,\alpha_2>0$ and $x_1,x_2\in B_X$. It follows that
  $\alpha_1x_1+\alpha_2x_2=\bigl(2-(\alpha_1+\alpha_2)\bigr)e$.  Thus,
  $\bigabs{2-(\alpha_1+\alpha_2)}\norm{e}=
  \norm{\alpha_1x_1+\alpha_2x_2}\leq\alpha_1+\alpha_2$.  If
  $\alpha_1+\alpha_2>2$ then
  $\bigl((\alpha_1+\alpha_2)-2\bigr)\norm{e}-(\alpha_1+\alpha_2)\leq
  0$, so that $\alpha_1+\alpha_2\leq\frac{2\norm{e}}{\norm{e}-1}$.  It
  follows that $\alpha_1\leq\alpha_1+\alpha_2\leq
  \max\bigl\{2,\frac{2\norm{e}}{\norm{e}-1}\bigr\}$ Finally, since
  $\norm{w}\leq\alpha_1(\norm{e}+1)+\norm{e}$, it follows that $W$ is
  bounded. Thus, $W$ is the unit ball of a norm, which is
  equivalent to the original norm of $X$. In the new norm 
  $e$ will be of norm one. Finally, $e$ dominates $W$ with respect to
  the order defined by $C$. Now apply Theorem~\ref{t:dom-cone}.
\end{proof}

\begin{remark}
  One can easily see that 
  Theorem~\ref{t:Krein} is equivalent to the original theorem of
  Krein. Indeed, if $T$ leaves invariant a non-empty open cone, then
  Theorem~\ref{t:Krein} states that $T^*$ has an eigenvector.
  Conversely, suppose that $T$ leaves invariant a cone with an
  interior point. Let $x$ be an interior point of the cone, then
  $f(x+Tx)\geq f(x)>0$ for every positive functional $f\neq 0$, so
  that $x+Tx$ is again an interior point of the cone. It follows that
  $I+T$ leaves invariant the interior of the cone, so that $(I+T)^*$
  has an eigenvector by Krein's theorem. This yields the existence of
  an eigenvector for $T^*$.
\end{remark}

Next, we discuss some applications of Theorems~\ref{t:dom-cone}
and~\ref{t:Krein}.

Recall that an element $e$ in a Banach lattice $E$ is called a
\term{strong unit} if for every positive $x\in E$ there exists a
natural number $n$ such that $x\leq ne$. It is known
(see~\cite[p.~188]{Aliprantis:85} for details) that a Banach lattice
with a strong unit is an AM-space with unit up to an equivalent norm.
But in an AM-space the unit dominates the unit ball. Therefore, 
Theorems~\ref{t:dom-cone} yields the following result.

\begin{corollary}
The adjoint of a positive
operator on a Banach lattice with a strong unit has a positive
eigenvector.
\end{corollary}

In particular, the adjoint of a positive operator on a
$C(\Omega)$ space, where $\Omega$ is a compact Hausdorff space, has a positive
eigenvector. A direct proof of this fact can also be found
in~\cite{Abramovich:98}.

The case of complex normed spaces can often be reduced to the real
case as follows. Suppose that $X_c$ is a complexification of a real
ordered normed space $X$, every element of $X_c$ can be written in the
form $x+iy$ for some $x,y\in X$.  If $T$ is a positive operator on
$X$, then its complexification $T_c\colon X_c\to X_c$ defined by
$T_c(x+iy)=Tx+iTy$ will be referred as a positive operator on $X_c$.
Notice that $T$ coincides with the restriction of $T_c$ to $X$.
Suppose that $T^*f=\lambda f$ for some $f\in X^*$ and
$\lambda\in\mathbb R$, then we can extend $f$ to a continuous linear
functional $f_c$ on $X_c$ via $f_c(x+iy)=f(x)+if(y)$. Then
$T_c^*f_c=\lambda f_c$. Indeed,
\begin{multline*}
  [T_c^*f_c](x+iy)=f_c\bigl(T_c(x+iy)\bigr)=f_c(Tx+iTy)=\\
  f(Tx)+if(Ty)=(T^*f)(x)+i(T^*f)(y)=\lambda f(x)+i\lambda f(y)=\lambda f_c(x+iy).
\end{multline*}
Thus, Theorems~\ref{t:dom-cone} and~\ref{t:Krein} are applicable to complex
normed spaces.

For example, we can apply our technique to positive operators on
$C^*$-algebras. A $C^*$-algebra $\iA$ can be viewed as
the complexification of the real Banach space $\iA_{\rm sa}$ of its
self-adjoint elements.  Recall that a self-adjoint element $a$ in
$\iA$ is positive if $\sigma(a)\subset\mathbb R_+$. If $\iA$ has unit $e$ and
$x$ is a self-adjoint element of $\iA$ such that $\norm{x}\leq 1$, then
the Spectral Mapping Theorem implies that $\sigma(e-x)\subseteq [0,2]$, hence
$x\leq e$. It follows that $e$ dominates the unit ball of $\iA_{\rm
 sa}$.  Theorem~\ref{t:dom-cone} immediately yields the following
result.

\begin{corollary} \lab{t:c-star}
  If $T$ is a positive operator on a unital $C^*$-algebra, then
  $T^*$ has a positive eigenvector.
\end{corollary}

Let $(e_j)_{j=1}^\infty$ denote the standard unit basis of $X=\ell_1$,
while $(e_i^*)_{i=1}^\infty$ stands for the dual basis of $X^*$. Recall
that every bounded operator $T$ on $\ell_1$ can be written as an
infinite matrix with entries $t_{ij}=\langle e_i^*,Te_j\rangle$.

\begin{theorem} \lab{t:l-1}
  Suppose that $T$ is a bounded operator on $\ell_1$ with matrix
  $(t_{ij})$, and suppose that there exists an index $k$ such that
  \begin{equation}
    \label{eq:tkk}
    t_{kk}\pm t_{kj}\geq\sum_{i \neq k}\abs{t_{ik} \pm t_{ij}}
  \end{equation}
  for each $j>1$. Then $T^*$ has a positive eigenvector.
\end{theorem}

\begin{proof}
  Without loss of generality $k=1$. Let $K$ be the cone spanned by
  $e_1+B_X$. It is easy to see that $K$ is spanned by the set $\{e_1\pm
  e_i\}_{i=2}^\infty$. We claim that
  $K=\{(x_i)\mid x_1\geq\sum_{i=2}^\infty\abs{x_i}\}$.  Indeed, it is easy to see
  that the later set is closed under addition and positive scalar
  multiplication, hence it is a cone. Furthermore, it contains
  $e_1\pm e_i$ for each $i\geq 2$, so that it contains $K$. Finally, if a
  non-zero sequence $(x_i)$ satisfies $x_1\geq\sum_{i=2}^\infty\abs{x_i}$ then
  $\frac{x}{x_1}-e_1\in B_X$, so that $(x_i)$ is contained in $K$.

  Clearly, $e_1$ dominates the unit ball of $X$ with respect to the
  order induced by $K$.
  The condition $t_{11}\pm t_{1j}\geq\sum_{i=2}^\infty\abs{t_{i1} \pm t_{ij}}$
  means that
  $$T(e_1\pm e_2)=Te_1\pm Te_j=(\text{the 1st column of }T) \pm (\text{the $j$-th column
    of }T)\in K$$
  for every $j\geq 2$, it follows that $T(K)\subseteq K$.
  Theorem~\ref{t:dom-cone} finishes the proof.
\end{proof}

\begin{example}
  Let $K$ be as in the preceding proof, and let $\mathcal C$ be
  the set of all operators on $\ell_1$ preserving $K$. Clearly, the
  adjoint of every operator in $\mathcal C$ has an eigenvector. By
  construction, $\mathcal C$ is itself a cone and a multiplicative
  semi-group in $\mathcal L(\ell_1)$. It is easy to see that $\mathcal
  C$ is closed in the strong operator topology (and, being a convex
  set, it is also closed in the weak operator topology). Finally, we
  claim that $\mathcal C$ has non-empty interior with respect to the
  norm topology of $\mathcal L(\ell_1)$. For example, put $S=(s_{ij})$
  such that $s_{ij}$ equals 1 if $i=j=1$ and 0 otherwise. We claim
  that $S$ is an interior point of $\mathcal C$. Indeed, suppose that
  $R=(r_{ij})$ such that $\norm{R}<\frac{1}{5}$, and let $T=S+R$. Show
  that $T\in\mathcal C$. Note that
  $\sum_{i=1}^\infty\abs{r_{ij}}=\norm{Re_j}<\frac{1}{5}$ for every $j\geq 1$.
  It follows that
  $t_{11}\pm t_{1j}=1+r_{11}\pm r_{1j}\geq
  1-\frac{1}{5}-\frac{1}{5}=\frac{3}{5}$
  for every $j>1$. On the other hand,
  $$\sum_{i=2}^\infty\abs{t_{i1} \pm t_{ij}}=\sum_{i=2}^\infty\abs{r_{i1} \pm r_{ij}}
    \leq\sum_{i=2}^\infty\abs{r_{i1}}+\sum_{i=2}^\infty\abs{r_{ij}}<\tfrac{2}{5}.$$
  Hence, $T$ satisfies (\ref{eq:tkk}) and, therefore, $T\in\mathcal C$.
\end{example}

\begin{corollary}
  If $T$ is an operator on a real Banach space leaving invariant a
  convex set whose interior is non-void and doesn't contain zero, then
  then $T^*$ has an eigenvector. Moreover, a commutative collection of
  such operators has a common eigenvector.
\end{corollary}

\begin{proof}
  Apply Theorem~\ref{t:Krein} to the cone generated by the invariant set.
\end{proof}

Krein's theorem gives a natural insight and provides a simple solution
to Exercise~VII.5.10 of \cite{Dunford:58}, even though at the first
glance the statement seems to have no connection to order structures.

\begin{proposition}[{\cite[Exercise VII.5.10]{Dunford:58}}] \label{p:DS}
  If $\norm{T}=1$ and $T$ has a non-zero fixed point, then $T^*$ has an
  eigenvector.
\end{proposition}

\begin{proof}
  Suppose that $\norm{T}=1$ and $Te=e$ for some $e$ of norm one.
  Then the set $e+B_X$ is invariant under $T$, so it the cone
  generated by this set. Clearly, this cone is proper and has non-void
  interior. Now apply Theorem~\ref{t:Krein}.
\end{proof}

This approach can be generalized as follows.

\begin{definition}
  If $X$ is an ordered normed space, we say that it has \term{monotone 
    norm} if $0\leq x\leq y$ implies $\norm{x}\leq\norm{y}$.
\end{definition}

\begin{theorem} \label{t:Te>e}
  Suppose that $T$ is a positive operator on an ordered normed space
  with monotone norm such that $\norm{T}=1$ and $Te\geq e$ for some
  $e>0$. Then $T^*$ has a positive eigenvector.
\end{theorem}

\begin{proof}
  Without loss of generality we can assume $\norm{e}=1$. Since the
  norm is monotone, we have $(e+X_+)\cap B_X^{\circ}=\varnothing$, where
  $B_X^{\circ}$ stands for the open unit ball of $X$.  Hence, $X_+\cap
  (B_X^{\circ}-e)=\varnothing$, so that the two sets can be separated by
  a positive functional $f$. Then $f$ is non-negative on $e-B_X$. Let
  $K$ be the closed cone generated by $X_+$ and $e-B_X$. Since $f$ is
  non-negative on $K$, it is, indeed, a proper cone.

  If $x\in B_X$ then $e-x\in K$, so that $e$ dominates $B_X$ in the
  (semi)order induced on $X$ by $K$. It is given that $T(X_+)\subseteq
  X_+\subseteq K$. Furthermore, if $x\in B_X$ then
  $Tx\in B_X$, and we have $T(e-x)=(Te-e)+(e-Tx)\in X_++(e-B_X)\subseteq
  K$, so that $T(e-B_X)\subseteq K$. It follows that $T(K)\subseteq
  K$. Now apply Theorem~\ref{t:dom-cone} to the order induced by $K$.
\end{proof}

Notice that the condition $\norm{T}=1$ in Theorem~\ref{t:Te>e} cannot be
dropped. Indeed, for any $\alpha>1$, let $T$ be $\alpha$ times the
left shift on $\ell_p$, $1\le p<\infty$, that is,
$T(x_1,x_2,x_3,\dots)=(\alpha x_2,\alpha x_3,\dots)$. Then
$\norm{T}=\alpha$ and $(1,\alpha^{-1},\alpha^{-2},\dots)$ is a fixed
point of $T$. Nevertheless, $T^*$ clearly has no eigenvectors.

It follows immediately that under the hypothesis of
Theorem~\ref{t:Te>e} the operator $T$ has an invariant subspace of
co-dimension one. In fact, there is a closed face of the positive cone
of $X$ which is invariant under $T$.  Recall that $E \subset X_+$ ($X_+$ is
the positive cone of $X$) is called a face of $X_+$ if $E$ is itself a
closed cone, and, for $x_1, x_2 \in X_+$, $x_1 + x_2 \in E$ implies
$x_1, x_2 \in E$.  One can easily see that a closed cone $E \subset X_+$ is
a face of $X_+$ iff it is hereditary, that is, $x \in E$ whenever $0 \leq
x \leq y$ and $y \in E$.

\begin{theorem} \label{face-cone}
  Suppose that $X$ is an ordered normed space with monotone norm and
  $T$ is a positive operator on a $X$ such that $\norm{T}=1$ and
  $Te>e$ for some $e>0$.  Then there exists a non-trivial closed face
  $E$ of the positive cone of $X$ which is invariant under $T$.
  Moreover, if $X$ is a Banach lattice then $E - E$ is closed
  non-trivial ideal in $X$, invariant under $T$.
\end{theorem}

\begin{proof}
  Without loss of generality $\norm{e}=1$. Let
  $$
  E=\{x\geq 0\mid 
  \lim_{\alpha \raw 0^+} (\norm{e+\alpha x} - 1)/\alpha = 0 \} .
  $$
  Note
  that if $x\in E$ and $0\leq y\leq x$ then $y\in E$. Note also that $E$ is non-trivial
  as the positive vector $Te-e\in E$ because, for
  $\alpha \in (0, 1)$,
  $$
  1 = \norm{e} \leq \norm{e+\alpha (Te-e)} \leq 
  \norm{(1-\alpha) e + \alpha Te} \leq
  (1-\alpha) \norm{e} +  \alpha\norm{Te} = 1.
  $$
  Furthermore, $E$ is $T$-invariant.
  Indeed, suppose $\alpha > 0$ and $x \in E$. Then
  $$
  \norm{e+\alpha T x} \leq \norm{T e+\alpha T x} \leq \norm{e+\alpha x} .
  $$
  Therefore,
  $$
  \lim_{\alpha \raw 0^+} (\norm{e+\alpha T x} - 1)/\alpha \leq
  \lim_{\alpha \raw 0^+} (\norm{e+\alpha x} - 1)/\alpha = 0.
  $$

  It is easy to see that $E$ is a cone.
  Indeed, if $x,y \in E$, then $c x \in E$ for $c > 0$, and
  $$
  \norm{e+\alpha(x+y)/2} - 1 \leq 
  \tfrac{1}{2}\bigl((\norm{e+\alpha x} - 1) + (\norm{e+\alpha y} - 1)\bigr) = o(\alpha) 
  $$
  as $\alpha$ approaches $0$.
  Thus, $x + y \in E$.

  To show that $E$ is closed, suppose $x_i$ is a sequence of positive elements
  in $E$, converging to $x$ in norm. We shall show that $x \in E$.
  Fix $\vr > 0$. It suffices to prove that, whenever $\alpha > 0$ is
  sufficiently small, the inequality
  $\norm{e + \alpha x} \leq 1 + \vr \alpha$ is satisfied.
  Find $i$ for which $\norm{x - x_i} < \vr/2$.
  There exists $\alpha_0$ such that
  $\norm{e + \alpha x_i} \leq 1 + \vr\alpha/2$ whenever $0<\alpha<\alpha_0$.
  Thus, for $\alpha \in (0, \alpha_0)$,
  $$
  \norm{e + \alpha x} \leq \norm{e + \alpha x_i} + \alpha \norm{x - x_i}
  \leq (1 + \vr \alpha/2) + \vr \alpha/2 = 1 + \vr \alpha .
  $$
  Finally, $e \notin E$, hence $E$ is a non-trivial face of the
  positive cone of $X$.

  Next, suppose that $X$ is a Banach lattice, and put $F = E-E$.
  Clearly, $F$ is an order ideal, that is, $F$ is a linear subspace
  such that $x\in F$ and $\abs{y}\leq\abs{x}$ imply $y\in F$.  Show that
  $F$ is closed.  Suppose $z \in \overline{F}$, and $(x_i)$, $(y_i)$
  are sequences in $E$ such that $\lim_i \norm{z - (x_i - y_i)} = 0$.
  Then $\lim_i \norm{z_+ - (x_i - y_i)_+} = 0$.  Let $a_i = (x_i -
  y_i)_+ \land z_+$. By the above, $\lim_i \norm{a_i - z_+} = 0$.  Note
  that
  $$
  0\leq a_i \leq (x_i - y_i)_+ \leq |x_i| + |y_i| \in E ,
  $$
  hence $a_i \in E$. But $E$ is closed, thus $z_+ \in E$.
  Similarly, $z_- \in E$, and therefore, $z \in F$.

  Finally we prove that $F$ is non-trivial.
  More precisely, we show that $e \notin F$.
  Indeed, suppose there exist $x,y\in E$ such that $\norm{e-(x-y)}\leq\frac{1}{3}$.
  Then $(x-y)_+\leq\abs{x}+\abs{y}\in E$, so $(x-y)_+\in E$. 
  Pick $\alpha > 0$ for which 
  $\norm{e + \alpha (x-y)_+} \leq 1 + \alpha/3$.
  Then
  \begin{multline*}
    1+\alpha=\norm{e+\alpha e}=
    \bignorm{e+\alpha(x-y)_++\alpha\bigl(e-(x-y)_+\bigr)}\\
    \leq\bignorm{e+\alpha(x-y)_+}+\alpha\bignorm{e-(x-y)_+}
    \leq 1+ \alpha/3 + \alpha\bignorm{e-(x-y)}=1+\tfrac{2\alpha}{3},
  \end{multline*}
  contradiction.
\end{proof}

Similar results hold for rearrangement invariant operator spaces,
arising from von Neumann algebras.  For the benefit of the reader, we
give a brief introduction into this natural non-commutative
generalization of Banach lattices.

Suppose $N$ is a von Neumann algebra on a Hilbert space $H$,
equipped with a faithful normal semifinite trace $\tau$.
Following~\cite{Nelson}, we say that a closed, densely defined linear operator
$x$ on $H$ is \term{affiliated with $N$} if $u^* x u = x$ for every unitary
$u \in N^\prime$ (the commutant of $N$).
An operator $x$ is called \term{$\tau$-measurable} if for every $\vr > 0$ there exists
a (self-adjoint) projection $p \in N$ such that $p(H) \subset D(x)$ and
$\tau(\one - p) < \vr$ ($\one$ is the identity in $N$).
The set of all $\tau$-measurable operators is denoted by $\tilde{N}$.

Following~\cite{Fack}, we introduce for $x \in \tilde{N}$
the \term{generalized eigenvalue function}
$\mu( \cdot , x) : [0, \infty) \to [0, \infty)$, defined by
$$
\mu(t,x) = \inf \{ s \geq 0  \mid
\tau( \chi_{(s, \infty)}(|x|) ) \leq t \} .
$$
Equivalently (see~\cite{Fack}), we have
$$
\mu(t,x) = \inf \{ \norm{xp}  \mid
p \in N \, a \, projection, \tau(\one - p) \leq t \} .
$$
Following~\cite{Dodds:93}, we call a linear manifold
$G \subset \tilde{N}$, equipped with the norm $\norm{\cdot}$,
a \term{(normed) rearrangement invariant operator space}
(\term{r.i.o.s.}, in short)
if whenever $x \in G$, $y \in \tilde{N}$,
and $\mu(t,y) \leq \mu(t,x)$ for every $t$, then $y \in G$ and
$\norm{y} \leq \norm{x}$.
$E$ is called \term{symmetric} if, in addition,
$\norm{y} \leq \norm{x}$ whenever
$$
\int_0^a \mu(t,y) dt \leq \int_0^a \mu(t,x) dt
$$
for every $a > 0$.

To underscore the connections between r.i.o.s. and Banach lattices,
consider the commutative case of $N = L_\infty(I)$,
where $I$ is an interval $(0, a)$ ($a \in (0, \infty]$).
By Proposition 2.a.8 of~\cite{LT:v2}, any r.i.o.s.~$G$ which satisfies
$$
L_1(I) \cap L_\infty(I) \subset E \subset L_1(I) + L_\infty(I)
\eqno(*)
$$
is symmetric.  We say that $G$ has the \term{Fatou property} if,
whenever $f \in G$, $(f_n)$ is a sequence of non-negative elements of
$G$, and $f_n(\omega) \nearrow f(\omega)$, then $\norm{f_n} \to
\norm{f}$.

Suppose $N$ is a von Neumann algebra with a normal faithful semifinite
trace $\tau$, and $G$ is as in the previous paragraph
(with $I = (0, \tau(\one))$).
Following \cite{Dodds:93}, we define the space
$G(N) = \{ x \in \tilde{N} \, | \, \mu( \cdot, x ) \in G \}$,
equipped with the norm $\norm{x}_{G(N)} = \norm{\mu( \cdot, x )}_G$.
If $G$ satisfies (*), then $N \cap N_* \subset G(N) \subset N + N_*$.
We identify $L_\infty(N)$ with $N$ itself, and $L_1(N)$ with $N_*$
(the predual of $N$).
If, in addition, $G$ has Fatou property, then $G(N)$ is norm closed
(see Proposition 1.7 and Corollary 2.4 of \cite{Dodds:93}).

If $G \subset N + N_*$ is a r.i.o.s., we denote by $G_+$ the set
of positive elements in $G$, i.e. $G \cap \tilde{N}_+$.
Then every self-adjoint element in $G$ can be represented as a difference
of two positive ones (see \cite{Dodds:89} and \cite{Dodds:93}).
Moreover, every element $x \in G$ can be written as
$x = x_1 - x_2 + i(x_3 - x_4)$, with $x_j \in G_+$.
Finally, the trace $\tau$ extends naturally to $(N + N_*)_+$
by setting $\tau(x) = \int_0^\infty \mu(t,x) dt$ for $x \geq 0$.

\begin{theorem} \label{noncommutative}
  Suppose $N$ is a von Neumann algebra with a faithful normal semifinite trace
  $\tau$, $G$ is a norm closed symmetric rearrangement invariant subspace of
  $\tilde{N}$ satisfying $N \cap N_* \subset G \subset N + N_*$, and
  $T : G \raw G$ is a positive contraction such that $Te > e$ for some
  positive $e \in G$.
  Then $T$ has an invariant non-trivial face $E$ of the positive cone of $G$
  Moreover, $\overline{E - E}$ is a non-trivial closed subspace of $G$,
  invariant under $T$.
\end{theorem}

To prove the theorem, we need to collect some facts related to
conditional expectations on von Neumann algebras.
Suppose $N$ is a von Neumann algebra equipped with a normal faithful
semifinite trace $\tau$, and $M$ is a Neumann subalgebra of $N$
s.t. the restriction of $\tau$ to $M$ is semifinite. Then 
(see Proposition~V.2.36 of \cite{Takesaki}) there
exists a positive contractive projection $\Phi$ from $N$ onto $M$ s.t.
$\Phi(a b c) = a \Phi(b) c$ and $\tau(\Phi(a b)) = \tau(a \Phi(b))$
whenever $a, c \in N_*$ and $b \in N$.
Moreover, it follows from the proof that, for any $x \in N \cap N_*$,
$\Phi(x) \in M \cap M_*$ and $\norm{\Phi(x)}_{M_*} \leq \norm{x}_{N_*}$.
Since $N \cap N_*$ (or $M \cap M_*$) is dense in $N_*$ (respectively, $M_*$),
$\Phi$ can be extended to a contraction from $N_*$ to $M_*$.
%%%   (we denote this extension by $\Phi$).
Thus, $\Phi$ can be thought of as an operator from $N + N_*$ to
$M + M_*$ respectively, which maps $N$ to $M$ and $N_*$ to $M_*$
contractively.

\begin{lemma} \label{cond-exp}
  Suppose $N$, $M$ and $\tau$ are as above, and $G$ is a symmetric r.i.o.s.
  with $N \cap N_* \subset G \subset N + N_*$.
  Then $\Phi$ maps $G$ into $G \cap \tilde{M}$, and
  $\norm{\Phi(x)}_G \leq \norm{x}_G$ for any $x \in G$.
\end{lemma}

\begin{proof}
  As noted above, $\Phi$ acts contractively from $N$ to $M$ and from $N_*$
  to $M_*$.
  For $x \in G$ we have, by Theorem 4.7 of \cite{Dodds:93},
  $\int_0^a \mu(t, \Phi(x)) d t$ $\leq \int_0^a \mu(t, x) d t$ for any $a > 0$.
  Thus, $\Phi(x) \in G$, and $\norm{\Phi(x)} \leq \norm{x}$.
\end{proof}

\begin{proof} [Proof of Theorem \ref{noncommutative}]
  %%%   We borrow some ideas from the proof of Theorem \ref{face-cone}.
  Suppose $e \in G_+$, $\norm{e} = 1$, and $T : G \raw G$ is a positive
  operator s.t. $Te > e$.
  Let 
  $$
  E=\{x\geq 0\mid 
  \lim_{\alpha \raw 0^+} (\norm{e+\alpha x} - 1)/\alpha = 0 \} .
  $$
  As in the proof of Theorem \ref{face-cone}, we can show that $E$ is
  a closed non-trivial face of $G_+$ ($E$ is non-empty, and $e \notin E$).
  Moreover, $E$ is invariant under $T$.
  Therefore the closed linear span of $E$ is invariant under $T$.
  It remains to show that $e$ does not belong to the closed linear span of $E$.
  It suffices to show that, whenever $x_1, x_2 \in E$, we have
  $\norm{e + x_1 - x_2} \geq 1/6$.

  First suppose that either $\tau(\one) < \infty$, or
  $\lim_{t \raw \infty} \mu(t, e) = 0$.
  Then there exists a commutative von Neumann algebra $M$ s.t.
  $e \in \tilde{M}$ and the restriction of $\tau$ to $M$ is
  semifinite.
  Indeed, if $\tau(\one) < \infty$, we can consider the von Neumann
  algebra generated by projections $\chi_{(a,\infty)}(e)$, where $a>0$.
  If $\lim_{t \raw \infty} \mu(t, e) = 0$, observe that
  $\tau(\chi_{(a,\infty)}(e)) < \infty$ for any $a > 0$, and let
  $p = \sup_{a>0} \chi_{(a,\infty)}(e)$.
  Use Zorn's lemma to find mutually orthogonal projections $(p_i) \in N$ s.t.
  $\tau(p_i) < \infty$ and $\sum_i p_i = \one - p$.
  Then let $M$ be the von Neumann algebra generated by projections
  $\chi_{(a,\infty)}(e)$ and $p_i$.
  Clearly $M$ satisfies our conditions.

  Let $\Phi$ be the conditional expectation from $N$ onto $M$.
  By Lemma~\ref{cond-exp}, $\Phi$ acts as a contraction from
  $G$ to $G_1 = G \cap \tilde{M}$.
  Then $G_1$ can be regarded as a Banach lattice.
  Let
  $$
  E_1=\{x \in G_1 \mid x \geq 0,
  \lim_{\alpha \raw 0^+} (\norm{e+\alpha x} - 1)/\alpha = 0 \} .
  $$
  As in the proof of Theorem~\ref{face-cone},
  $\norm{e + x - y} \geq 1/3$ whenever $x, y \in E_1$.
  However, $\Phi(E) \subset E_1$, and therefore
  $$
  \norm{e + x - y} \geq \norm{e + \Phi(x) - \Phi(y)} \geq \frac{1}{3}
  $$
  whenever $x, y \in E$.

  The case of $a = \lim_{t \raw \infty} \mu(t,e) > 0$ is more complicated.
  Note that $\norm{a \one}_G \leq \norm{e} = 1$, hence
  $\norm{x}_G \leq \norm{x}_N \norm{\one}_G \leq \norm{x}_N/a$
  for any $x \in N$.
  Let $k = \lceil 6/a \rceil$, $p_i = \chi_{[ia/k, (i+1)a/k)}(e)$ for
  $0 \leq i \leq k-1$, $p_k = \chi_{[(k-1)a/k, a]}(e)$, and
  $e_1 =
  \chi_{(a,\infty)}(e) e + \sum_{i=1}^k \frac{i}{k} a p_i$.
  Then $e \geq e_1$, $e - e_1 \in N$, and
  $$
  \norm{e - e_1}_G \leq \norm{e - e_1}_N/a \leq 1/6 .
  $$
  By definition, $\mu(t, e) = \mu(t, e_1)$ for any $t$.
  Moreover, a projection $p_i$ can be represented as
  $p_i = \sum_j q_{ij}$, where projections $q_{ij}$ are mutually orthogonal
  and $\tau(q_{ij}) < \infty$.
  Note also that $\tau(\chi_{(b,\infty)}(e)) < \infty$ whenever $b > a$,
  and $\chi_{(a,\infty)}(e) = \sup_{b > a} \chi_{(b,\infty)}(e)$.
  Consider the (commutative) von Neumann algebra $M$, generated by
  projections $q_{ij}$ and $\chi_{(b,\infty)}(e)$ ($b > a$).
  Then $e_1 \in G_1 = G \cap \tilde{M}$. Let
  $$
  E_1 = \{x \in G_1 \mid x \geq 0,
  \lim_{\alpha \raw 0^+} (\norm{e_1 + \alpha x} - 1)/\alpha = 0 \} .
  $$
  As above, we show that $\norm{e_1 + x - y} \geq 1/3$
  if $x, y \in E_1$.
  However, $\Phi(e) \geq e_1$ (since $\Phi$ is positive), and
  therefore, 
  $\norm{e_1 + \Phi(x)} \leq \norm{e + x}$ for any $x \in G$.
  Thus, $\Phi(E) \in E_1$ and, for any $x, y \in E$,
  we have
  \begin{displaymath}
  \norm{e + x - y} \geq \norm{\Phi(e) + \Phi(x) - \Phi(y)}
      \geq
  \norm{e_1 + \Phi(x) - \Phi(y)} - \norm{e - e_1} \geq
  \frac{1}{3} - \frac{1}{6}.
  \end{displaymath}
  The proof is complete.
\end{proof}

Thanks are due to Y.A.~Abramovich and C.D.~Aliprantis for their
interest in the work and useful suggestions. The authors would like to
thank the organizers of the workshop in Linear Analysis in College
Station in 2001 where much of the work on this paper was done. The
authors also thank the Department of Mathematics of the University of
Texas at Austin for its hospitality and support.


\begin{thebibliography}{DDdP89}

\bibitem[AA02]{AA:02}
Y.~A. Abramovich and C.~D. Aliprantis
\newblock {\em An Invitation to Operator Theory}.
\newblock Graduate Studies in Mathematics,
Amer.\ Math.\ Soc., Providence, RI, 2002.

\bibitem[AAB98]{Abramovich:98}
Y.~A. Abramovich, C.~D. Aliprantis, and O.~Burkinshaw.
\newblock The invariant subspace problem: Some recent advances.
\newblock {\em Rend. Inst. Mat. Univ. Trieste}, XXIX Supplemento:3--79, 1998.

\bibitem[AAB92]{Abramovich:92}
Y.~A. Abramovich, C.~D. Aliprantis, and O.~Burkinshaw.
\newblock Positive operators on {K}re\u\i n spaces.
\newblock {\em Acta Appl. Math.}, 27(1-2):1--22, 1992.

\bibitem[AB85]{Aliprantis:85}
C.~D. Aliprantis and O.~Burkinshaw.
\newblock {\em Positive operators}.
\newblock Academic Press Inc., Orlando, Fla., 1985.

\bibitem[DDdP89]{Dodds:89}
P.~Dodds, T.~Dodds, and B.~dePagter.
\newblock Noncommutative Banach function spaces.
\newblock {\em Math. Z.}, 201:583--597, 1989.

\bibitem[DDdP93]{Dodds:93}
P.~Dodds, T.~Dodds, and B.~dePagter.
\newblock Noncommutative K\"othe duality.
\newblock {\em Trans. Amer. Math. Soc.}, 339:717--750, 1993.

\bibitem[DS58]{Dunford:58}
N.~Dunford and J.T.~Schwartz.
\newblock {\em Linear {O}perators {I}}.
\newblock Interscience Publishers, Inc., New York, 1958.

\bibitem[FK86]{Fack}
T.~Fack and H.~Kosaki.
\newblock Generalized $s$-numbers of $\tau$-measurable operators.
\newblock {\em Pacific J. Math.}, 123:269--300, 1986.

\bibitem[KR48]{Krein:48}
M.~G. Kre{\u{\i}}n and M.~A. Rutman.
\newblock Linear operators leaving invariant a cone in a {B}anach space.
\newblock {\em Uspehi Matem. Nauk (N. S.)}, 3(1(23)):3--95, 1948.
\newblock English translation:
%
%\bibitem[KR50]{Krein:50}
%M.~G. Kre{\u{\i}}n and M.~A. Rutman.
%\newblock Linear operators leaving invariant a cone in a {B}anach space.
\newblock {\em Amer. Math. Soc. Translation}, 1950(26):128, 1950.

\bibitem[LT79]{LT:v2}
Y.~Lindenstrauss and L.~Tzafriri.
\newblock {\em Classical Banach spaces II}.
\newblock Springer-Verlag, Berlin-New York, 1979.

\bibitem[N74]{Nelson}
E.~Nelson.
\newblock Notes on non-commutative integration.
\newblock {\em J. Functional Analysis}, 15:103--116, 1974.

\bibitem[S99]{Schaefer}
Schaefer, H. H. with Wolff, M. P.
\newblock {\em Topological vector spaces}.
\newblock Springer-Verlag, New York, 1999.

\bibitem[T79]{Takesaki}
M.~Takesaki.
\newblock {\em Theory of operator algebras I}.
\newblock Springer-Verlag, New York-Heidelberg, 1979.

\end{thebibliography}
\end{document}